\documentclass[12pt]{amsart}
\usepackage{fullpage,url}
\usepackage[all]{xy} 

\DeclareFontEncoding{OT2}{}{} 
\newcommand{\textcyr}[1]{%
 {\fontencoding{OT2}\fontfamily{wncyr}\fontseries{m}\fontshape{n}\selectfont #1}}
\newcommand{\Sha}{{\mbox{\textcyr{Sh}}}}

\usepackage{color}

\def\max{\mathop{\rm max}}

\def\cha{\mathop{\rm char}}

\def\pmb#1{\setbox0=\hbox{#1}%
 \kern-.025em\copy0\kern-\wd0
 \kern.05em\copy0\kern-\wd0
 \kern-.025em\raise.0433em\box0 }

\def\Z{{\bf Z}}
 
\def\Q{{\bf Q}}
\def\A{{\bf A}}

\def\F{{\bf F}}

\def\calO{{\mathcal O}}
\def\isom{\cong}

\def\Gal{\mathop{\rm Gal}\nolimits}
\def\Pic{\mathop{\rm Pic}\nolimits}

\DeclareMathOperator{\Div}{Div}
\DeclareMathOperator{\divv}{div}
\DeclareMathOperator{\Sel}{Sel}
\DeclareMathOperator{\Aut}{Aut}

\DeclareMathOperator{\BM}{BM}
\DeclareMathOperator{\Tate}{Tate}
\DeclareMathOperator{\Br}{Br}
\DeclareMathOperator{\ev}{ev}

\newcommand{\Abar}{{\overline{A}}}
\newcommand{\Fbar}{{\overline{F}}}
\newcommand{\Kbar}{{\overline{K}}}
\newcommand{\Xbar}{{\overline{X}}}
\newcommand{\Adual}{{A^\vee}}
\newcommand{\Selhat}{{\widehat{\Sel}}}
\newcommand{\Directsum}{\bigoplus}
\newcommand{\Union}{\bigcup}
\newcommand{\Intersection}{\bigcap}

\newcommand{\injects}{\hookrightarrow}
\newcommand{\intersect}{\cap}
\newcommand{\tensor}{\otimes}
\newcommand{\directsum}{\oplus}
\newcommand{\Zhat}{{\hat{\Z}}}
\newcommand{\calA}{{\mathcal A}}
\newcommand{\calB}{{\mathcal B}}
\newcommand{\calT}{{\mathcal T}}

\newcommand{\calY}{{\mathcal Y}}

\newcommand{\G}{{\mathbf G}}
\newcommand{\kbar}{{\overline{k}}}
\newcommand{\OO}{{\mathcal O}}

\newcommand{\tors}{{\operatorname{tors}}}
\newcommand{\all}{{\operatorname{all}}}

\newtheorem{theorem}{Theorem}
   
\newtheorem{lemma}{Lemma}[section]
\newtheorem{corollary}[lemma]{Corollary}
\newtheorem{proposition}[lemma]{Proposition}

\theoremstyle{definition}

\newtheorem{conjecture}[theorem]{Conjecture}

\theoremstyle{remark}
\newtheorem{remark}[lemma]{Remark}

\usepackage[
	backref,
	pdfauthor={Bjorn Poonen and Jose Felipe Voloch}, 
	pdftitle={The Brauer-Manin obstruction for subvarieties of abelian varieties over function fields},
]{hyperref}
\usepackage[alphabetic,backrefs,lite]{amsrefs} 

\begin{document}
\title[Brauer-Manin obstruction]{The Brauer-Manin obstruction for subvarieties of abelian varieties over function fields}
\author{Bjorn Poonen}
\address{Department of Mathematics, University of California, 
	Berkeley, CA 94720-3840, USA}
\email{poonen@math.berkeley.edu}
\urladdr{http://math.berkeley.edu/\~{}poonen}

\author{Jos\'e Felipe Voloch}
\address{Department of Mathematics, University of Texas,
	Austin, TX 78712, USA}
\email{voloch@math.utexas.edu}
\urladdr{http://www.ma.utexas.edu/\~{}voloch}

\date{December 7, 2006}

\begin{abstract}
We prove that for a large class of subvarieties of abelian varieties
over global function fields, the Brauer-Manin condition on adelic points
cuts out exactly the rational points.
This result is obtained from more general results concerning
the intersection of the adelic points of a subvariety with the 
adelic closure of the group of rational points of the abelian variety.
\end{abstract}

\maketitle

\section{Introduction}

The notation in this section remains in force throughout the paper,
except in Section~\ref{S:Mordell-Lang},
and in Section~\ref{S:Brauer-Manin} where we allow also the possibility
that $K$ is a number field.

Let $k$ be a field.
Let $K$ be a finitely generated extension of $k$ of transcendence degree $1$.
We assume that $k$ is relatively algebraically closed in $K$,
since the content of our theorems will be unaffected by this assumption.
Let $\Kbar$ be an algebraic closure of $K$.
We will use this notation consistently for an algebraic closure,
and we will choose algebraic closures compatibly whenever possible:
thus $\kbar$ is the algebraic closure of $k$ in $\Kbar$.
Let $K^s$ be the separable closure of $K$ in $\Kbar$.
Let $\Omega_\all$ be the set of all nontrivial valuations 
of $K$ that are trivial on $k$.
Let $\Omega$ be a cofinite subset of $\Omega_\all$;
if $k$ is finite, we may weaken the cofiniteness hypothesis to assume only
that $\Omega \subseteq \Omega_\all$ has Dirichlet density $1$.
For each $v \in \Omega$, let $K_v$ be the completion of $K$ at $v$,
and let $\F_v$ be the residue field.
Equip $K_v$ with the $v$-adic topology.
Define the ring of ad\`eles $\A$ as the restricted direct product
$\prod_{v \in \Omega} (K_v,\calO_v)$
of the $K_v$ with respect to their valuation subrings $\calO_v$.
Then $\A$ is a topological ring, in which $\prod_{v \in \Omega} \calO_v$ 
is open and has the product topology.

If $A$ is an abelian variety over $K$,
then $A(K)$ embeds diagonally into $A(\A) \simeq \prod_v A(K_v)$.
Define the adelic topology on $A(K)$ as the topology induced from $A(\A)$.
For any fixed $v$
define the $v$-adic topology on $A(K)$ as the topology induced from $A(K_v)$.
Let $\overline{A(K)}$ be the closure of $A(K)$ in $A(\A)$.

For any extension of fields $F' \supset F$ and any $F$-variety $X$,
let $X_{F'}$ be the base extension of $X$ to $F'$.
Call a $K$-variety $X$ {\em constant}
if $X \isom Y_K$ for some $k$-variety $Y$,
and call $X$ {\em isotrivial} if $X_\Kbar \isom Y_\Kbar$ 
for some variety $Y$ defined over $\kbar$.

{}From now on, $X$ is a closed $K$-subscheme of $A$.
Call $X$ {\em coset-free} if 
$X_{\Kbar}$ does not contain a translate of a positive-dimensional
subvariety of $A_{\Kbar}$.

When $k$ is finite and $\Omega=\Omega_\all$,
the intersection $X(\A) \cap {\overline {A(K)}} \subset A(\A)$
is closely related to 
the Brauer-Manin obstruction to the Hasse principle for $X/K$:
see Section~\ref{S:Brauer-Manin}.
For curves over number fields,
V. Scharaschkin and A. Skorobogatov independently raised the question 
of whether the Brauer-Manin obstruction 
is the only obstruction to the Hasse principle,
and proved that this is so when the Jacobian has finite Mordell-Weil group
and finite Shafarevich-Tate group.
The connection with the adelic intersection is stated explicitly
in \cite{Scharaschkin-thesis},
and is based on global duality statements 
originating in the work of Cassels: see Remark~\ref{R:history}.
See also \cites{Skorobogatov2001,Flynn2004,Poonen-heuristic2005preprint,Stoll2006preprint}: the last of these in particular contains many conjectures 
and theorems relating descent information,
the method of Chabauty and Coleman, 
the Brauer-Manin obstruction,
and Grothendieck's section conjecture.

In this paper we answer (most cases of) a generalization of
the function field analogue 
of a question raised for curves over number fields
in \cite{Scharaschkin-thesis},
concerning whether the Brauer-Manin condition cuts out exactly
the rational points: see Theorem~\ref{T:Selmer}.
This question is still wide open in the number field case.
Along the way, we prove results about adelic intersections
along the lines of the ``adelic Mordell-Lang conjecture''
suggested in \cite{Stoll2006preprint}*{Question~3.12}.
Again, these are open in the number field case.
In particular, we prove the following theorems.

\begin{theorem}
\label{T:A}
If $\cha k = 0$, then $X(K)= X(\A) \cap {\overline {A(K)}}$.
\end{theorem}

\begin{theorem}
\label{T:B}
Suppose that $\cha k = p > 0$, 
that $A_{\Kbar}$ has no nonzero isotrivial quotient,
and that $A(K^s)[p^\infty]$ is finite.
Suppose that $X$ is coset-free.
Then $X(K)=X(\A) \cap {\overline{A(K)}}$.
\end{theorem}

\begin{remark}
The proposition in \cite{Voloch1995} states that 
in the ``general case'' in which $A$ is ordinary and 
the Kodaira-Spencer class of $A/K$ has maximal rank,
we have $A(K^s)[p^\infty]=0$.
\end{remark}

\begin{conjecture}
\label{C:conjecture}
For any closed $K$-subscheme $X$ of any $A$, we have
$\overline{X(K)}=X(\A) \cap {\overline{A(K)}}$,
where $\overline{X(K)}$ is the closure of $X(K)$ in $X(\A)$.
\end{conjecture}

\begin{remark}
If $A_{\Kbar}$ has no nonzero isotrivial quotient
and $X$ is coset-free,
then $X(K)$ is finite~\cite{Hrushovski1996}*{Theorem~1.1}
so $\overline{X(K)}=X(K)$.
So Conjecture~\ref{C:conjecture} predicts in particular that
the hypothesis on $A(K^s)[p^\infty]$ in Theorem~\ref{T:B}
is unnecessary.
\end{remark}

\begin{remark}
Here is an example to show that the statement
$X(K)=X(\A) \cap {\overline{A(K)}}$
can fail for a constant curve in its Jacobian.
Let $X$ be a curve of genus $\ge 2$ over a finite field $k$.
Choose a divisor of degree $1$ on $X$ 
to embed $X$ in its Jacobian $A$.
Let $F\colon A \to A$ be the $k$-Frobenius map.
Let $K$ be the function field of $X$.
Let $P \in X(K)$ be the point corresponding to the identity map $X \to X$.
Let $P_v \in X(\F_v)$ be the reduction of $P$ at $v$.

For each $v$, the Teichm\"uller map $\F_v \to K_v$ identifies $\F_v$ 
with a subfield of $K_v$.
Any $Q \in A(K_v)$ can be written as $Q=Q_0+Q_1$
with $Q_0 \in A(\F_v)$ and $Q_1$ in the kernel of 
the reduction map $A(K_v) \to A(\F_v)$;
then $\lim_{m \to \infty} F^m(Q_1)=0$,
so $\lim_{n \to \infty} F^{n!}(Q)=\lim_{n \to \infty} F^{n!}(Q_0) = Q_0$.
In particular, taking $Q=P$,
we find that $\left(F^{n!}(P) \right)_{n \ge 1}$ converges in $A(\A)$
to the point $(P_v) \in X(\A) = \prod_v X(K_v)$,
where we have identified $P_v$ with its image under
the Teichm\"uller map $X(\F_v) \injects X(K_v)$.
If $(P_v)$ were in $X(K)$, then in $X(K_v)$ we would have
$P_v \in X(\F_v) \cap X(K) = X(k)$,
which contradicts the definition of $P_v$ 
if $v$ is a place of degree greater than $1$ over $k$.
Thus $(P_v)$ is in $X(\A) \cap {\overline{A(K)}}$ but not in $X(K)$.
\end{remark}

In the final section of this paper, we restrict to
the case of a global function field,
and extend Theorem~\ref{T:B} to prove (under mild hypotheses)
that for a subvariety of an abelian variety, 
the Brauer-Manin condition cuts out exactly the rational points:
see Section~\ref{S:Brauer-Manin} for the definitions of
$X(\A)^{\Br}$ and $\Selhat$.
Our result is as follows.

\begin{theorem}
\label{T:Selmer}
Suppose that $K$ is a global function field of characteristic $p$,
that $A_{\Kbar}$ has no nonzero isotrivial quotient,
and that $A(K^s)[p^\infty]$ is finite.
Suppose that $X$ is coset-free.
Then $X(K) = X(\A)^{\Br} = X(\A) \intersect \Selhat$.
\end{theorem}

To our knowledge, Theorem~\ref{T:Selmer} is the first 
result giving a wide class of varieties of general type 
such that the Brauer-Manin condition cuts out exactly the rational points.

\section{Characteristic 0}

Throughout this section, we assume $\cha k=0$.
In this case, results follow rather easily.

\begin{proposition}
\label{P:v-adic=discrete}
For any $v$, the $v$-adic topology on $A(K)$ equals
the discrete topology.
\end{proposition}

\begin{proof}
The question is isogeny-invariant, so we reduce to the case 
where $A$ is simple.
Let $A(\F_v)$ denote the group of $\F_v$-points on the N\'eron model
of $A$ over $\calO_v$.
Let $A^1(K_v)$ be the kernel of the reduction map
$A(K_v) \to A(\F_v)$.
The Lang-N\'eron theorem \cite{Lang-Neron1959}*{Theorem~1}
implies that either $A$ is constant and $A(K)/A(k)$ is finitely generated,
or $A$ is non-constant and $A(K)$ itself is finitely generated.
In either case, the subgroup $A^1(K) := A(K) \cap A^1(K_v)$ 
is finitely generated.
Now $A^1(K_v)$ has a descending filtration by open subgroups 
in which the quotients of consecutive terms are torsion-free
(this is where we use $\cha k=0$), 
so the induced filtration on the finitely generated group $A^1(K)$
has only finitely many nonzero quotients.
Thus $A^1(K)$ is discrete.
Since $A^1(K_v)$ is open in $A(K_v)$,
the group $A(K)$ is discrete in $A(K_v)$.
\end{proof}

\begin{remark}
The literature contains results close to Proposition~\ref{P:v-adic=discrete}.
It is mentioned in the third subsection of the introduction 
to~\cite{Manin1963} for elliptic curves with non-constant $j$-invariant,
and it appears for abelian varieties with $K/k$-trace zero 
as the main theorem of \cite{Buium-Voloch1993}.
\nocite{Manin1989}
\end{remark}

\begin{corollary}
\label{C:adelic=discrete}
The adelic topology on $A(K)$ equals
the discrete topology.
\end{corollary}

\begin{proof}
The adelic topology is at least as strong as 
(i.e., has at least as many open sets as)
the $v$-adic topology for any $v$.
\end{proof}

We can improve the result by imposing conditions
in only the residue fields $\F_v$ 
instead of the completions $K_v$, that is, ``flat'' instead of ``deep'' 
information in the sense of \cite{Flynn2004}.
In fact, we have:

\begin{proposition}
There exist $v,v' \in \Omega$ of good reduction for $A$
such that $A(K) \to A(\F_v) \times A(\F_{v'})$ is injective.
\end{proposition}

\begin{proof}
Let $B$ be the $K/k$-trace of $A$.
Pick any $v \in \Omega$ of good reduction.
The kernel $H$ of $A(K) \to A(\F_v)$ meets $B(k)$ trivially.
By Silverman's specialization theorem 
\cite{LangFundamentals}*{Chapter~12, Theorem~2.3},
there exists $v' \in \Omega$ such that $H$ injects under
reduction modulo $v'$.
\end{proof}

\begin{proof}[Proof of Theorem~\ref{T:A}]
By Corollary~\ref{C:adelic=discrete}, 
$X(\A) \cap {\overline {A(K)}} = X(\A) \cap  A(K)= X(K)$.
\end{proof}

\section{Characteristic $p$}
\label{S:char p}

Throughout this section, $\cha k=p$.

\subsection{Abelian varieties}
\label{S:abelian varieties}

The following is a function field analogue 
of \cite{Serre1971}*{Th\'eor\`eme 3}.
It is also essentially a weak version of the main theorem 
of \cite{Milne1972b}.

\begin{proposition}
\label{P:prime-to-p}
The adelic topology on $A(K)$ 
is at least as strong as 
the topology induced by all subgroups of finite prime-to-$p$ index.
\end{proposition}

\begin{proof}
As in the proof of Proposition~\ref{P:v-adic=discrete},
the Lang-N\'eron theorem implies that $A(\A)$
has an open subgroup 
whose intersection with $A(K)$ is finitely generated.
It suffices to study the topology induced on that
finitely generated subgroup,
so we may reduce to the case 
in which $k$ is finitely generated over a finite field $\F_q$.
In particular, now $A(K)$ is finitely generated.
It suffices to prove that for each prime $\ell \ne p$, 
the adelic topology on $A(K)$
is at least as strong as the topology induced 
by subgroups of $\ell$-power index.

We use the terminology of \cite{Serre1964}, but with $\ell$ in place of $p$.
Let $T_\ell:=\varprojlim_n A(K^s)[\ell^n]$ 
be the $\ell$-adic Tate module of $A$, 
and let $V_\ell:=T_\ell \tensor_{\Z_\ell} \Q_\ell$.
Let $G_\ell$ be the image of $\Gal(K^s/K) \to \Aut_{\Q_\ell} V_\ell$,
and let ${\mathfrak g}_\ell$ be the corresponding Lie algebra.
Define
\[
	H^1_S(K,A[\ell^n]):= \ker\left( H^1(K,A[\ell^n]) \to \prod_{v \in \Omega} H^1(K_v,A[\ell^n]) \right).
\]
By \cite{Serre1971}*{Remarque~2, p.~734} applied to $K$,
the Lie algebra cohomology group $H^1({\mathfrak g}_\ell,V_\ell)$ vanishes.
The proof of \cite{Serre1964}*{Th\'eor\`eme~3} 
deduces from this that $\varprojlim_n H^1_S(K,A[\ell^n])=0$:
the one place in that proof where it is used that $K$ is a number field
is in the proof of \cite{Serre1964}*{Proposition~8},
which uses that a cyclic subgroup of a finite quotient of $\Gal(K^s/K)$
is the image of a decomposition group $D_v$ for some $v \in \Omega$,
but it is still true in our setting (and it suffices for the argument)
that such a cyclic subgroup
is {\em contained} in the image of a decomposition group;
this is because the decomposition group contains Frobenius elements
attached to closed points of an $\F_q$-variety with function field $K$
(or rather its finite extension)
and there is a version of the Chebotarev density theorem for finitely
generated fields (cf.~\cite{Serre1965}*{Theorem 7}).
The analogue of \cite{Serre1964}*{Th\'eor\`eme~2} 
(whose proof proceeds as in the
number field case, again making use of the analogue of 
\cite{Serre1964}*{Proposition~8} just discussed)
says that the (just-proved) vanishing of $\varprojlim_n H^1_S(K,A[\ell^n])$
implies the desired conclusion: that for every subgroup $H$ of 
$\ell$-power index in $A(K)$,
there is an open subgroup $U$ of $A(\A)$
such that $A(K) \intersect U \subseteq H$.
\end{proof}

\begin{lemma}
\label{L:K_v separable}
If $\alpha \in K_v$ is algebraic over $K$,
then $\alpha$ is separable over $K$.
\end{lemma}

\begin{proof}
Replacing $K$ by its relative separable closure in $L:=K(\alpha)$,
we may assume that $L$ is purely inseparable over $K$.
Then the valuation $v$ on $K$ admits a unique extension $w$ to $L$,
and the inclusion of completions $K_v \to L_w$ is an isomorphism.
By \cite{SerreLocalFields}*{I.\S4, Proposition~10}
(hypothesis (F) there holds for localizations of 
finitely generated algebras over a field),
we have an ``$n=\sum e_i f_i$'' result,
which in our case says $[L:K] = [L_w:K_v] = 1$.
So $\alpha \in K$.
\end{proof}

\begin{lemma}
\label{L:Hausdorff}
For any $n \in \Z_{\ge 1}$,
the quotient $A(K_v)/n A(K_v)$ is Hausdorff.
\end{lemma}

\begin{proof}
Equivalently, we must show that $n A(K_v)$ is closed in $A(K_v)$.
Suppose $(P_i)$ is a sequence in $n A(K_v)$ that converges to $P \in A(K_v)$.
Write $P_i=n Q_i$ with $Q_i \in A(K_v)$.
Then $n(Q_i-Q_{i+1}) \to 0$ as $i \to \infty$.

Let $\OO_v$ be the valuation ring of $K_v$,
and let $\calA$ over be the N\'eron model of $A$ over $\OO_v$.
Applying \cite{Greenberg1966}*{Corollary~1} to $\calA[n]$
shows that for any sequence $(R_i)$ in $A(K_v)$ with $n R_i \to 0$,
the distance of $R_i$ to the nearest point of $A(K_v)[n]$ tends to $0$.

Thus by induction on $i$ we may adjust each $Q_i$ by a point in $A(K_v)[n]$
so that $Q_i-Q_{i+1} \to 0$ as $i \to \infty$.
Since $A(K_v)$ is complete, $(Q_i)$ converges to some $Q \in A(K_v)$,
and $n Q = P$.
Thus $n A(K_v)$ is closed.
\end{proof}

The following is a slight generalization of \cite{Voloch1995}*{Lemma~2},
with a more elementary proof.

\begin{proposition}
\label{P:p-power}
If $A(K^s)[p^\infty]$ is finite,
then for any $v$, the $v$-adic topology on $A(K)$ 
is at least as strong as 
the topology induced by all subgroups of finite $p$-power index.
\end{proposition}

\begin{proof}
For convenience choose algebraic closures $\overline{K},\overline{K}_v$ 
of $K,K_v$ such that $K^s \subseteq \overline{K} \subseteq \overline{K}_v$.
As in the proof of Proposition~\ref{P:v-adic=discrete},
there is an open subgroup $U$ of $A(K_v)$ such that 
$B:=A(K) \intersect U$ is finitely generated.
It suffices to show that for every $e \in \Z_{\ge 0}$, 
there exists an open subgroup $V$ of $A(K_v)$ such that 
$B \intersect V \subseteq p^e A(K)$.

Choose $m$ such that $p^m A(K^s)[p^\infty] = 0$.
Let $M=e+m$.
Then $B/p^M B$ is finite.
By Lemma~\ref{L:Hausdorff}, $A(K_v)/p^M A(K_v)$ is Hausdorff,
so the image of $B/p^M B$ in $A(K_v)/p^M A(K_v)$ is discrete.
Hence there is an open subgroup $V$ of $A(K_v)$
such that $B \intersect V = \ker(B \to A(K_v)/p^M A(K_v))$.

Suppose $b \in B \intersect V$.
Then $b=p^M c$ for some
$c \in A(K_v) \intersect A(\overline{K})$.
By Lemma~\ref{L:K_v separable}, we obtain $c \in A(K^s)$.
If $\sigma \in \Gal(K^s/K)$, then ${}^\sigma c - c \in A(K^s)[p^M]$
is killed by $p^m$.
Thus $p^m c \in A(K)$.
So $b = p^e p^m c \in p^e A(K)$.
Hence $B \intersect V \subseteq p^e A(K)$.
\end{proof}

If $A(K^s)[p^\infty]$ is finite,
we can simply combine Propositions \ref{P:prime-to-p} and~\ref{P:p-power}
to obtain the following, but we will prove it even without an assumption
on $A(K^s)[p^\infty]$.

\begin{proposition}
\label{P:finite index}
The adelic topology on $A(K)$ 
is at least as strong as 
the topology induced by all subgroups of finite index.
\end{proposition}

\begin{proof}
In \cite{Milne1972b}, the result is proved for the case 
where $k$ is finitely generated.  
(There instead of the adelic topology as we have defined it,
he uses the topology coming from the closed points of a finite-type
$\Z$-scheme with function field $K$, but the adelic topology is stronger,
so his result contains what we want.)
The general case can be reduced to the case 
where $A(K)$ is finitely generated,
and hence to the case where $k$ is finitely generated,
by using the same argument used at the beginning of the proofs
of Propositions \ref{P:prime-to-p} and~\ref{P:p-power}.
\end{proof}

\begin{lemma}
\label{L:topological group}
Suppose $A(K)$ is finitely generated.
Then $\overline{A(K)}_\tors = A(K)_\tors$.
\end{lemma}

\begin{proof}
Let 
\[
	T := \ker\left(A(K) \to \prod_{v \in \Omega} \frac{A(\F_v)}{A(\F_v)_\tors} \right),
\]
where $A(\F_v)$ is the group of $\F_v$-points on the N\'eron model of $A$.
Since $A(K)$ is finitely generated and the groups $A(\F_v)/A(\F_v)_\tors$
are torsion-free, there is a finite subset $S \subset \Omega$ such that 
$T = A(K) \intersect U$ for the open subgroup
\[
	U := \ker\left(A(\A) \to \prod_{v \in S} 
				\frac{A(\F_v)}{A(\F_v)_\tors} \right)
\]
of $A(\A)$.
The finitely generated group $A(K)/T$ is contained in 
the torsion-free group $\prod_{v \in S} \frac{A(\F_v)}{A(\F_v)_\tors}$,
so $A(K)/T$ is free,
and we have $A(K) \isom T \directsum F$
as topological groups, where $F$ is a discrete free abelian group
of finite rank.

We claim that the topology of $T$ is that induced
by the subgroups $nT$ for $n \ge 1$.
For $n \ge 1$, the subgroup $nT$ is open in $T$
by Proposition~\ref{P:finite index}.
If $t \in T$, then some positive multiple of $t$ is in the kernel of
$A(K_v) \to A(\F_v)$, and then $p$-power multiples of this multiple
tend to $0$.
Applying this to a finite set of generators of $T$,
we see that any open neighborhood of $0$ in $T$ contains $nT$
for some $n \in \Z_{>0}$.

It follows that $\overline{T} \isom T \tensor \Zhat$.
Now
\[
	\overline{A(K)}_\tors = (\overline{T} \directsum F)_\tors 
	= \overline{T}_\tors \isom (T \tensor \Zhat)_\tors \isom T_\tors.
\]
\end{proof}

The following proposition is a function field analogue of
\cite{Stoll2006preprint}*{Proposition~3.6}.
Our proof must be somewhat different, however, since \cite{Stoll2006preprint}
made use of strong ``image of Galois'' theorems 
whose function field analogues have recently been 
disproved \cite{Zarhin2006preprint}.

\begin{proposition}
\label{P:subscheme}
Suppose that $A(K^s)[p^\infty]$ is finite.
Let $Z$ be a finite $K$-subscheme of $A$.
Then $Z(\A) \cap \overline{A(K)} = Z(K)$.
\end{proposition}

\begin{proof}
As in \cite{Stoll2006preprint}, 
we may replace $K$ by a finite extension 
to assume that $Z$ consists of a finite set of $K$-points of $A$.
A point $P \in \overline{A(K)}$ 
is represented by a sequence $(P_n)_{n \ge 1}$ in $A(K)$
such that for every $v$, the limit $\lim_{n \to \infty} P_n$
exists in $A(K_v)$.
If in addition $P \in Z(\A)$, then there is a point $Q_v \in Z(K)$ 
whose image in $Z(K_v)$ equals $\lim_{n \to \infty} P_n \in A(K_v)$.
The $P_n-Q_v$ are eventually contained in the kernel of
$A(K) \to A(\F_v)$, which is finitely generated,
so there are finitely generated subfields $k_0 \subseteq k$, 
$K_0 \subseteq K$ with $K_0/k_0$ a function field 
such that all the $P_n$ and the points of $Z(K)$ are in $A(K_0)$.
By Proposition~\ref{P:p-power},
the sequence $(P_n-Q_v)_{n \ge 1}$ is eventually divisible in $A(K_0)$ by
an arbitrarily high power of $p$.
For any other $v' \in \Omega$,
then the same is true of $(P_n-Q_{v'})_{n \ge 1}$.
Then $Q_{v'}-Q_v \in A(K_0)$ is divisible by every power of $p$.
Since $A(K_0)$ is finitely generated, 
$Q_{v'}-Q_v$ is a torsion point in $A(K_0)$.
This holds for every $v' \in \Omega$, and $A(K_0)_\tors$ is finite,
so $R:=P-Q_v \in \overline{A(K_0)}$ is a torsion point in $\overline{A(K_0)}$.
Lemma~\ref{L:topological group} applied to $K_0$ yields $R \in A(K_0)_\tors$.
Thus $P = R+Q_v \in A(K)$.
So $P \in Z(\A) \cap A(K) = Z(K)$.
\end{proof}

\begin{lemma}
\label{L:closure mod p^e}
Fix $v \in \Omega$.
Let $\Gamma_v$ be the closure of $A(K)$ in $A(K_v)$.
Then for every $e \in \Z_{\ge 0}$, 
the map $A(K)/p^e A(K) \to \Gamma_v/p^e \Gamma_v$ is surjective.
\end{lemma}

\begin{proof}
The theory of formal groups implies that $A(K_v)$ 
has an open subgroup $A^\circ(K_v)$
that is a torsion-free topological $\Z_p$-module,
and we may choose $A^\circ(K_v)$ so that $A^\circ(K):=A(K) \cap A^\circ(K_v)$
is finitely generated.
The group $\Gamma_v^\circ:=\Gamma_v \cap A^\circ(K_v)$
is the closure of $A^\circ(K)$, 
so there is an isomorphism of topological groups
$\Gamma_v^\circ \isom \Z_p^{\oplus m}$ 
for some $m \in \Z_{\ge 0}$.
In particular, for any $e \in \Z_{\ge 0}$,
the group $p^e \Gamma_v^\circ$ is open in $\Gamma_v^\circ$,
which is open in $\Gamma_v$.
So the larger group $p^e \Gamma_v$ also is open in $\Gamma_v$.
But the image of $A(K)$ in the discrete group $\Gamma_v/p^e \Gamma_v$ is dense,
so the map $A(K)/p^e A(K) \to \Gamma_v/p^e \Gamma_v$ is surjective.
\end{proof}

\subsection{A uniform Mordell-Lang conjecture}
\label{S:Mordell-Lang}

We thank Zo\'e Chatzidakis, Fran\c{c}oise Delon, and Tom Scanlon
for many of the ideas used in this section.

\begin{lemma}
\label{L:infinite orbit}
Let $\calB$ be a $p$-basis for a field $K$ of characteristic $p$.
Let $L$ be an extension of $K$ such that $\calB$ is also a $p$-basis for $L$.
Suppose that $c$ is an element of $L$ that is not algebraic over $K$.
Then there exists a separably closed extension $F$ of $L$
such that $\calB$ is a $p$-basis of $F$ and the $\Aut(F/K)$-orbit of $c$
is infinite.
\end{lemma}

\begin{proof}
Fix a transcendence basis $T$ for $L/K$.
Let $\Omega$ be an algebraically closed extension of $K$
such that the transcendence basis of $\Omega/K$ is identified
with the set $\Z \times T$.
Identify $L$ with a subfield of $\Omega$ 
in such a way that each $t \in T$
maps to the transcendence basis element for $\Omega/K$
labelled by $(0,t) \in \Z \times T$.
The map of sets $\Z \times T \to \Z \times T$ mapping $(i,t)$ to $(i+1,t)$
extends to an automorphism $\sigma \in \Aut(\Omega/K)$.
For $i \in \Z$, let $L_i = \sigma^i(L)$.
Let $L_\infty$ be the compositum of the $L_i$ in $\Omega$.
Then $\sigma(L_\infty)=L_\infty$.
Let $F$ be the separable closure of $L_\infty$ in $\Omega$.
So $\sigma(F)=F$.
The $\sigma$-orbit of $c$ is infinite,
since $L_i \intersect L_j$ is algebraic over $K$ whenever
$i \ne j$.

The $p$-basis hypothesis implies that $L$ is separable over $K$
(see \cite{Delon1998} for the definitions of 
separable, $p$-basis, $p$-free, etc.)
Applying $\sigma^i$ shows that $L_i$ is separable over $K$.
Moreover, the $L_i$ are algebraically disjoint over $K$,
so their compositum $L_\infty$ is separable over $K$,
by the last corollary in \cite{BourbakiAlgebra4-7}*{V.\S 16.7}
and Proposition~3(b) in \cite{BourbakiAlgebra4-7}*{V.\S 15.2}.
Thus $\calB$ is $p$-free in $L_\infty$.

The $p$-basis hypothesis implies also that $L=L^p(\calB)$.
Thus $L_i=L_i^p(\calB)$ for all $i \in \Z$,
and $L_\infty = L_\infty^p(\calB)$.

Combining the previous two paragraphs shows that $\calB$ is a $p$-basis 
for $L_\infty$, and hence also for $F$.
\end{proof}

\begin{lemma}
\label{L:totally nonisotrivial}
Let $K \subseteq F$ be a separable extension
such that the field $K^{p^\infty}:=\Intersection_{n \ge 1} K^{p^n}$ 
is algebraically closed.
Let $A$ be an abelian variety over $K$ such that 
no nonzero quotient of $A_\Kbar$ is the base extension of
an abelian variety over $K^{p^\infty}$.
Then no nonzero quotient of $A_\Fbar$ is the base extension of
an abelian variety over $F^{p^\infty}$.
\end{lemma}

\begin{proof}
Suppose not.
Thus there exists a nonzero abelian variety $B$ over $F^{p^\infty}$
and a surjective homomorphism $\phi\colon A_\Fbar \to B_\Fbar$.
Choose a finitely generated extension $F_0$ of $K^{p^\infty}$
over which $B$ is defined.
Then $\phi$ is definable over a finite separable extension $F_1$ of $KF_0$,
since any homomorphism $A \to B$ between abelian varieties
defined over a field $L$ is definable over $L(A[\ell^\infty],B[\ell^\infty])$
for any prime $\ell$.
Since $K^{p^\infty}$ is algebraically closed,
we may choose a place $F_0 \dashrightarrow K^{p^\infty}$ 
extending the identity on $K^{p^\infty}$,
such that $B$ has good reduction at this place.
By \cite{Delon1998}*{Fact~1.4}, $K$ and $F^{p^\infty}$ are linearly disjoint
over $K^{p^\infty}$,
so the place $F_0 \dashrightarrow K^{p^\infty}$ 
extends to a place $KF_0 \dashrightarrow K$
that is the identity on $K$,
and then to a place $F_1 \dashrightarrow K_1$ for some finite extension 
$K_1$ of $K$.
Reduction of $\phi\colon A_{F_1} \to B_{F_1}$
yields a homomorphism $A_{K_1} \to B_{K_1}$,
and since the place restricted to $F_0$ has values in $K^{p^\infty}$,
the abelian variety $B_{K_1}$ is the base extension of an 
abelian variety over $K^{p^\infty}$.
This contradicts the hypothesis on $A$.
\end{proof}

\begin{lemma}[a version of the Mordell-Lang conjecture]
\label{L:Hrushovski}
Let $F$ be a separably closed field of characteristic $p$.
Suppose that $A$ is an abelian variety over $F$
such that no nonzero quotient of $A_{\overline{F}}$ is the
base extension of an abelian variety over $F^{p^\infty}$.
Suppose that $X$ is a coset-free closed $F$-subvariety of $A$.
Then there exists $e \in \Z_{\ge 1}$ such that for every $a \in A(F)$,
the intersection $X(F) \intersect (a + p^e A(F))$ is finite.
\end{lemma}

\begin{proof}
This is a special case of \cite{Hrushovski1996}*{Lemma~6.2}.
\end{proof}

Using Lemma~\ref{L:Hrushovski}
and the compactness theorem in model theory,
we can deduce a version of the 
Mordell-Lang conjecture that is more uniform as we vary the ground field:

\begin{proposition}[more uniform version of the Mordell-Lang conjecture]
\label{P:independent of F}
Let $k$ be an algebraically closed field of characteristic $p$.
Let $K$ be a finitely generated extension of $k$.
Fix a (finite) $p$-basis $\calB$ of $K$.
Suppose that $A$ is an abelian variety over $K$
such that no nonzero quotient of $A_{\Kbar}$ is 
the base extension of an abelian variety over $k$.
Suppose that $X$ is a coset-free closed $K$-subscheme of $A$.
Let $\pi\colon A \times X \to A$ be the first projection.
Let $\calY$ be the set of $K$-subvarieties $Y \subseteq A \times X$
such that $\pi|_Y \colon Y \to A$ has finite fibers.
Then there exists $e \in \Z_{\ge 1}$ and $Y \in \calY$
such that for every field extension $F \supseteq K$ 
having $\calB$ as $p$-basis,
if $a \in A(F)$ and $c \in X(F) \intersect (a + p^e A(F))$, 
then $(a,c) \in Y(F)$.
\end{proposition}

\begin{proof}
Consider the language of fields augmented 
by a constant symbol $\alpha_\kappa$ for each element $\kappa \in K$ 
and by additional finite tuples of constant symbols $a$ and $c$
(to represent coordinates of points on $A$ and $X$, respectively).
We construct a theory $\calT$ in this language.
Start with the field axioms,
and the arithmetic sentences involving the $\alpha_\kappa$
that hold in $K$.
Include a sentence saying that 
$\calB$ is a $p$-basis for the universe field $F$.
Include sentences saying that $a,c \in A(F)$
(for some fixed representation of $A$ as a definable set).
For each $e \in \Z_{\ge 1}$, include a sentence saying that
\[
	c \in X(F) \intersect (a + p^e A(F)).
\]
For each $Y \in \calY$, include a sentence saying that $(a,c) \notin Y(F)$.

Suppose that the conclusion of Proposition~\ref{P:independent of F} fails.
Then we claim that $\calT$ is consistent.
To prove this, we show that every finite subset $\calT_0$ of $\calT$
has a model.
Given $\calT_0$,
let $e_{\max}$ be the maximum $e$ that occurs in the sentences,
and let $Y_{\max} \in \calY$ be the union of the $Y$'s that occur.
The negation of Proposition~\ref{P:independent of F}
implies that there exists a field extension
$F \supseteq K$ having $\calB$ as a $p$-basis,
together with $a \in A(F)$ and $c \in X(F) \intersect (a + p^{e_{\max}} A(F))$,
such that $(a,c) \notin Y_{\max}(F)$.
This $(F,a,c)$ is a model for $\calT_0$.

By the compactness theorem, there is a model $(F,a,c)$ for $\calT$.
So $F$ is a field extension of $K$ having $\calB$ as a $p$-basis,
and
\[
	c \in X(F) \intersect (a + p^e A(F)).
\]
for all $e \in \Z_{\ge 1}$
and for every $Y \in \calY$ we have $(a,c) \notin Y(F)$.
Then $K(c)$ is not algebraic over $K(a)$, because an algebraic relation
could be used to define a $Y \in \calY$ containing $(a,c)$.
By Lemma~\ref{L:infinite orbit} applied to $F/K(a)$,
it is possible to enlarge $F$ to assume that $F$ is separably closed
and the $\Aut(F/K(a))$-orbit of $c$ is infinite 
while $\calB$ is still a $p$-basis for $F$.
This infinite orbit is contained in the set
$X(F) \intersect (a + p^e A(F))$,
so the latter is infinite too, for every $e \in \Z_{\ge 1}$.
By Lemma~\ref{L:totally nonisotrivial} applied to $K \subseteq F$,
the abelian variety $A_{\overline{F}}$ 
has no nonzero quotient defined over $F^{p^\infty}$.
Also, the fact that $X$ is coset-free implies that $X_F$ is coset-free
(if $X_{\Kbar}$ acquires a coset over some extension of $\Kbar$,
it acquires one over some finitely generated extension of $\Kbar$,
and then a standard specialization argument shows
that it contains a coset already over $\Kbar$).
Thus $X_F \subseteq A_F$ is a counterexample to Lemma~\ref{L:Hrushovski}.
\end{proof}

\subsection{Subvarieties of abelian varieties}

We now return to our situation in which $K$ is the function field
of a curve over $k$.

\begin{lemma}
\label{L:exists Z}
Suppose that $A_{\Kbar}$ has no nonzero isotrivial quotient,
and that $X$ is coset-free.
Then there exists a finite $K$-subscheme $Z$ of $X$
such that $X(\A) \cap \overline{A(K)} \subseteq Z(\A)$.
\end{lemma}

\begin{proof}
By replacing $K$ with its compositum with $\kbar$ in $\Kbar$,
we may reduce to the case in which $k$ is algebraically closed.

We have $[K:K^p]=p$.
Choose $\alpha \in K-K^p$.
Then $\{\alpha\}$ is a $p$-basis for $K$.
For any $v$, the field $K_v$ is generated by $K$ and $K_v^p$,
so $[K_v:K_v^p] \le p$.
Moreover, Lemma~\ref{L:K_v separable} implies $\alpha \notin K_v^p$,
so $\{\alpha\}$ is a $p$-basis for $K_v$.

Let $e$ and $Y$ be as in Proposition~\ref{P:independent of F}.
For each $a \in A(K)$, let $Y_a$ be the fiber of $\pi|_Y \colon Y \to A$
above $a$, viewed as a finite subscheme of $X$.

Since $A_{\Kbar}$ has no nonzero isotrivial quotient, 
$A(K)$ is finitely generated.
Choose a (finite) set of representatives $\calA \subseteq A(K)$ 
for $A(K)/p^e A(K)$.
Let $Z = \Union_{a \in \calA} Y_a$, so $Z$ is a finite subscheme of $X$.
The conclusion of Proposition~\ref{P:independent of F}
applied to $F=K_v$ says that
\[
	X(K_v) \intersect (a + p^e A(K_v)) \subseteq Y_a(K_v)
\]
for each $a \in \calA$.
Since $\Gamma_v \subseteq A(K_v)$, we obtain
\[
	X(K_v) \intersect (a + p^e \Gamma_v) \subseteq Y_a(K_v).
\]
Taking the union over $a \in \calA$ yields
\[
	X(K_v) \intersect \Gamma_v \subseteq Z(K_v),
\]
since by Lemma~\ref{L:closure mod p^e}, 
the $a \in \calA$ represent all classes in $\Gamma_v/p^e \Gamma_v$.
This holds for all $v$, so $X(\A) \cap \overline{A(K)} \subseteq Z(\A)$.
\end{proof}

\subsection{End of proof of Theorem~\ref{T:B}}

We now prove Theorem~\ref{T:B}.
Let $Z$ be as in Lemma~\ref{L:exists Z}.
By Lemma~\ref{L:exists Z} and then Proposition~\ref{P:subscheme},
\[
	X(\A) \intersect \overline{A(K)} 
	\subseteq Z(\A) \cap \overline{A(K)} 
	= Z(K) 
	\subseteq X(K).
\]
The opposite inclusion $X(K) \subseteq X(\A) \intersect \overline{A(K)}$
is trivial.
This completes the proof.

\begin{remark}
Our proof of Theorem~\ref{T:B} required both the flat information
(from residue fields) in Proposition~\ref{P:prime-to-p}
and the deep information (from each $K_v$) in Proposition~\ref{P:p-power}.
But it seems possible that the flat information suffices,
even for Conjecture~\ref{C:conjecture}.
Specifically, if $\Omega$ 
is a cofinite subset of the set of places of good reduction for $X/K$,
then it is possible that 
$\overline{X(K)}=X\left(\prod_{v \in \Omega} \F_v \right) \cap {\overline{A(K)}}$
always holds, where the closures are now taken in
$A\left(\prod_{v \in \Omega} \F_v \right)$.
\end{remark}

\section{The Brauer-Manin obstruction}
\label{S:Brauer-Manin}

In this section we assume that $K$ is either 
a global function field (i.e., as before, but with $k$ finite)
or a number field.
We assume also that $\Omega$ is the set of all nontrivial places of $K$,
and define the ring of ad\`eles $\A$ in the usual way.

The purpose of this section is to relate the adelic intersections we
have been considering to the Brauer-Manin obstruction.
This relationship was discovered by Scharaschkin \cite{Scharaschkin-thesis} 
in the number field case.
We will give a proof similar in spirit to his, 
and verify that it works in the function field case.

For convenience, if $X$ is any topological space,
let $X_\bullet$ denote the set of connected components,
and equip $X_\bullet$ with the quotient topology:
this will be needed to avoid annoyances with the
archimedean places in the number field case.

Let $A$ be an abelian variety over $K$.
By the adelic topology on $A(K)$,
we mean the topology induced by $A(\A)_\bullet$;
this agrees with our previous terminology in the function field case.

\begin{proposition}
\label{P:finite index for global fields}
The adelic topology on $A(K)$ equals the topology induced by subgroups
of finite index.
\end{proposition}

\begin{proof}
Proposition~\ref{P:finite index} gives one inclusion.
The other follows from the fact that 
$A(\A)_\bullet=\prod_{v \in \Omega} A(K_v)_\bullet$
is a profinite group.
\end{proof}

\subsection{Global duality for abelian varieties}

All cohomology below is fppf cohomology.
As usual, for any $n \in \Z_{\ge 1}$, define
\begin{align*}
	\Sel^n &:= \ker\left( H^1(K,A[n]) \to \prod_{v \in \Omega} H^1(K_v,A) \right) \\
	\Sha &:=\ker\left( H^1(K,A) \to \prod_{v \in \Omega} H^1(K_v,A) \right).
\end{align*}
Following the notation of \cite{Stoll2006preprint}*{\S2}, define
\begin{align*}
	\widehat{A(K)} &:= \varprojlim \frac{A(K)}{nA(K)} \\
	\Selhat &:=\varprojlim \Sel^n \\
	T\Sha &:=\varprojlim \Sha[n],
\end{align*}
where each inverse limit is over positive integers $n$ ordered by divisibility.
Taking inverse limits of the usual descent sequence yields
\begin{equation}
\label{E:descent}
	0 \to \widehat{A(K)} \to \Selhat \to T\Sha \to 0.
\end{equation}

For any abelian profinite group $G$ and any prime $\ell$, 
let $G^{(\ell)}$ be the maximal pro-$\ell$ quotient of $G$,
so $G \isom \prod_\ell G^{(\ell)}$.

\begin{lemma}
\label{L:finitely generated}
For each prime $\ell$,
each of
$\widehat{A(K)}^{(\ell)}$, $\Selhat^{(\ell)}$, $T\Sha^{(\ell)}$
is a finitely generated $\Z_\ell$-module.
\end{lemma}

\begin{proof}
For $\widehat{A(K)}^{(\ell)}$ 
it holds simply because $A(K)$ is finitely generated as a $\Z$-module.
For $T\Sha^{(\ell)}$, it follows from the finiteness of $\Sha[\ell]$,
which is proved in \cite{Milne1970b}.
Now the result for $\Selhat^{(\ell)}$ follows from \eqref{E:descent}.
\end{proof}

By definition, $\Sel^n$ maps to
\[
	\ker\left( \prod_{v \in \Omega} H^1(K_v,A[n]) \to \prod_{v \in \Omega} H^1(K_v,A) \right) = \prod_{v \in \Omega} \frac{A(K_v)}{n A(K_v)},
\]
so as in \cite{Stoll2006preprint}*{\S2} we obtain a map
\[
	\Selhat \to A(\A)_\bullet.
\]

For each $v$ we have a pairing
\[
	A(K_v)_\bullet \times H^1(K_v,\Adual) \to H^2(K_v,\G_m) \isom \Q/\Z,
\]
and in fact Tate local duality holds: each of the first two groups
is identified with the Pontryagin dual
of the other \cite{MilneADT}*{Theorem~III.7.8}.
An element of $H^1(K,\Adual)$ maps to $0$ in $H^1(K_v,\Adual)$
for all but finitely many $v$, so summing the Tate local duality pairing
over $v$ defines a pairing
\[
	A(\A)_\bullet \times H^1(K,\Adual) \to \Q/\Z,
\]
or equivalently a homomorphism
\[
	A(\A)_\bullet \stackrel{\Tate}\longrightarrow H^1(K,\Adual)^D,
\]
where the superscript ${}^D$ denotes Pontryagin dual.

The following global duality statement
is a version of the ``Cassels dual exact sequence'' in which finiteness
of $\Sha$ is not assumed:
\begin{proposition}[Cassels dual exact sequence]
\label{P:Cassels dual exact sequence}
The sequence
\[
	0 \longrightarrow \Selhat \longrightarrow A(\A)_\bullet \stackrel{\Tate}\longrightarrow H^1(K,\Adual)^D
\]
is exact.
\end{proposition}

\begin{proof}
This is part of the main theorem of \cite{GonzalezAviles-Tan-preprint2006}.
\end{proof}

\begin{remark}
\label{R:history}
Proposition~\ref{P:Cassels dual exact sequence}
has a long history.
For an elliptic curve $A$ over a number field,
Proposition~\ref{P:Cassels dual exact sequence}
was proved by Cassels \cites{Cassels1962,Cassels1964}
under the assumption that $\Sha(A)$ is finite, 
using a result from \cite{Serre1964} to obtain exactness 
on the left.
This, except for the exactness on the left, 
was generalized by Tate to abelian varieties
over number fields with finite Shafarevich-Tate group,
at around the same time, though the first published proof 
appeared much later \cite{MilneADT}*{Theorem~I.6.13}.
The latter reference also proved the global function field analogue
except for the $p$-part in characteristic $p$.
The $\Selhat$ version (not requiring finiteness of $\Sha(A)$),
except for exactness on the left and the $p$-part, 
is implicit in the right three terms
of the exact sequence in the middle of page 104 of \cite{MilneADT}
if there we make the substitution
\[
	T \Sha \isom \Selhat / \widehat{A(K)}
\]
from \eqref{E:descent}.
Exactness on the left for the prime-to-$p$ part 
is implicit in \cites{Serre1964,Serre1971}
and is mentioned explicitly in \cite{MilneADT}*{Corollary~6.23}.
The full statement was proved in \cite{GonzalezAviles-Tan-preprint2006}.
\end{remark}

\subsection{Review of the Brauer-Manin obstruction}

Let $X$ be a smooth projective geometrically integral variety over $K$.
Let $\Xbar := X \times_K K^s$.
Define $\Br X := H^2(X,\G_m)$, where here we may use \'etale cohomology
since it gives the same result as fppf cohomology for smooth quasi-projective
commutative group schemes \cite{MilneEtale}*{Theorem~III.3.9}.
For simplicity, let $\frac{\Br X}{\Br K}$ 
denote the cokernel of $\Br K \to \Br X$ even if this map is not injective.
The Hochschild-Serre spectral sequence in \'etale cohomology gives
an exact sequence
\[
	\Br K \to \ker\left(\Br X \to \Br \Xbar\right) \to H^1(K,\Pic \Xbar) \to H^3(K,\G_m),
\]
and the last term is $0$ for any global field: 
see \cite{MilneADT}*{Corollary~4.21}, for instance.
{}From this we extract an injection
\[
	H^1(K,\Pic \Xbar) \injects \frac{\Br X}{\Br K};
\]
if $\Br \Xbar=0$, then this is an isomorphism.
Composing with the map induced by $\Pic^0 \Xbar \injects \Pic \Xbar$, 
we obtain
\begin{equation}
\label{E:H^1 of Pic^0}
	H^1(K,\Pic^0 \Xbar) \to \frac{\Br X}{\Br K}.
\end{equation}

Each $A \in \Br X$ induces evaluation maps $\ev_A$ fitting in a commutative
diagram
\[
\xymatrix{
& X(K) \ar[r] \ar[d]_{\ev_A} & X(\A)_\bullet \ar[d]_{\ev_A} \ar@{..>}[rd] \\
0 \ar[r] & \Br K \ar[r] & \displaystyle \Directsum_{v \in \Omega}^{\phantom{v \in \Omega}}  \Br(K_v) \ar[r] & \Q/\Z \ar[r] & 0. \\
}
\]
The diagonal arrows give rise to a pairing
\[
	X(\A)_\bullet \times \frac{\Br X}{\Br K} \to \Q/\Z
\]
which induces a map
\[
	X(\A)_\bullet \stackrel{\BM}\longrightarrow \left(\frac{\Br X}{\Br K}\right)^D.
\]
The kernel (inverse image of $0$) is called the {\em Brauer set},
and is denoted $X(\A)_\bullet^{\Br}$:
it contains $X(K)$
and hence also its closure $\overline{X(K)}$ in $X(\A)_\bullet$.
One says that the Brauer-Manin obstruction to the Hasse principle
is the only one for $X$ if the implication
\[
	X(\A)_\bullet^{\Br} \ne \emptyset \implies X(K) \ne \emptyset
\]
holds.
(This is equivalent to the usual notion without the $\bullet$.)

\subsection{The Brauer set for abelian varieties}

\begin{theorem}
\label{T:Brauer-Manin for A}
We have
\[
	\widehat{A(K)} = \overline{A(K)} \subseteq A(\A)_\bullet^{\Br} \subseteq \Selhat.
\]
in $A(\A)$.
\end{theorem}

\begin{proof}
Proposition~\ref{P:finite index for global fields}
shows that $A(K) \to A(\A)_\bullet$ induces an isomorphism
$\widehat{A(K)} \to \overline{A(K)}$,
and the first inclusion is automatic from the previous subsection.

So it remains to prove the second inclusion.
Taking the dual of \eqref{E:H^1 of Pic^0} for $A$ yields the vertical arrow in
\[
\xymatrix{
0 \ar[r] & \Selhat \ar[r] & A(\A)_\bullet \ar[r]^-{\Tate} \ar[rd]_-{\BM} & H^1(K,\Pic^0 \Abar)^D \\
&&& \displaystyle \left( \frac{\Br A}{\Br K} \right)^D \ar[u] \\
}
\]
Applying \cite{Manin1971}*{Proposition~8} to each $K_v$
shows that the triangle commutes
(the hypothesis there that the ground field be perfect is not used).
The diagram now gives $\ker(\BM) \subseteq \ker(\Tate)$,
or equivalently $A(\A)_\bullet^{\Br} \subseteq \Selhat$.
\end{proof}

\begin{remark}
If $\Sha$ is finite (or more generally, if its maximal divisible part
$\Sha_{\divv}$ is $0$), then $T\Sha=0$
and the injection $\widehat{A(K)} \to \Selhat$
is an isomorphism;
in this case, the conclusion of Theorem~\ref{T:Brauer-Manin for A}
simplifies to
\[
	\widehat{A(K)} = \overline{A(K)} = A(\A)_\bullet^{\Br}.
\]
Closely related results in the number field case 
can be found in \cite{Wang1996}.
\end{remark}

\subsection{The Brauer set for subvarieties of abelian varieties}

\begin{proposition}
\label{P:BM for subvarieties}
Let $X$ be a subvariety of $A$, and let $i \colon X \to A$ be the inclusion.
Then $i$ induces a map $X(\A)_\bullet \to A(\A)_\bullet$
(possibly non-injective if there are archimedean places).
\begin{enumerate}
\item[(a)]
We have $i(X(\A)^{\Br}_\bullet) \subseteq \Selhat$.
\item[(b)]
Suppose that
\begin{itemize}
\item $X$ is a smooth projective geometrically integral curve,
\item $A$ is the Jacobian of $X$,
\item $i$ is an Albanese map (that is, $i$ sends a point $P$ 
to the class of $P-D$,
where $D \in \Div \Abar$ is a fixed divisor of degree $1$ whose class
is $\Gal(K^s/K)$-invariant), and
\item $\Sha$ is finite (or at least $\Sha_{\divv}=0$).
\end{itemize}
Then $X(\A)^{\Br}_\bullet = i^{-1}(\overline{A(K)})$.
\end{enumerate}
\end{proposition}

\begin{proof}
(a) 
This follows immediately from 
Theorem~\ref{T:Brauer-Manin for A},
since $i$ maps $X(\A)^{\Br}_\bullet$ to $A(\A)^{\Br}_\bullet$.

(b) 
We use the diagram
\[
\xymatrix{
&&& \displaystyle \left( \frac{\Br X}{\Br K} \right)^D \ar[d] \ar@/^6pc/[dddd] \\
&&& H^1(K,\Pic \Xbar)^D \ar[d] \\
&& X(\A)_\bullet \ar[ruu]^-{\BM_X} \ar[d]_i & H^1(K,\Pic^0 \Xbar)^D \ar[d] \\
0 \ar[r] & \overline{A(K)} \ar[r] & A(\A)_\bullet \ar[r]^-{\Tate} \ar[rd]_-{\BM_A} & H^1(K,\Pic^0 \Abar)^D \\
&&& \displaystyle \left( \frac{\Br A}{\Br K} \right)^D \ar[u] \\
}
\]
in which the horizontal sequence is the exact sequence of
Proposition~\ref{P:Cassels dual exact sequence}
in which we have used the finiteness of $\Sha$ to
replace $\Selhat$ with $\overline{A(K)}$.
The lower triangle commutes: as mentioned earlier, 
this follows from a result of Manin.
The ``pentagon'' at the far right also commutes,
since the homomorphism \eqref{E:H^1 of Pic^0} is functorial in $X$.
Thus the whole diagram commutes.

We next claim that the three downward vertical arrows 
at the right are isomorphisms.
The first is an isomorphism because 
$\Br \Xbar=0$ \cite{Grothendieck-Brauer3}*{Corollaire 5.8}.
The second is an isomorphism 
because $\Pic \Xbar \isom \Pic^0 \Xbar \directsum \Z$
(where the $\Z$ is generated by the class of $D$)
and $H^1(K,\Z)=0$.
The third is an isomorphism because an Albanese map $i$ induces
an isomorphism $i^*\colon \Pic^0 \Abar \to \Pic^0 \Xbar$.

The commutativity of the upper left ``hexagon''
now implies
\[
	X(\A)^{\Br}_\bullet := \ker\left(\BM_X\right)
	= \ker(\Tate \circ i) = i^{-1}(\ker(\Tate)) = i^{-1}(\overline{A(K)}).
\]
\end{proof}

\begin{remark}
Part~(b) was originally proved in the number field case
in \cite{Scharaschkin-thesis} using a slightly different proof.
For yet another proof of this case, 
see \cite{Stoll2006preprint}*{Corollary~7.4}.
\end{remark}

\begin{remark}
If $K$ is a global function field, the conclusion of (b) can be written as
$X(\A)^{\Br} = X(\A) \intersect \overline{A(K)}$ in $A(\A)$.
\end{remark}

\section{Intersections with $\Selhat$}

{}From now on, we assume that $K$ is a global function field.

The proof of Theorem~\ref{T:Selmer} will follow that of Theorem~\ref{T:B},
with $\Selhat$ playing the role of $\overline{A(K)}$.
We begin by proving $\Selhat$-versions of several of the lemmas and 
propositions of Section~\ref{S:abelian varieties}.
The following is an analogue of Lemma~\ref{L:topological group}.

\begin{lemma}
\label{L:Selmer torsion}
The maps $A(K)_\tors \to \widehat{A(K)}_\tors \to \Selhat_\tors$
are isomorphisms.
\end{lemma}

\begin{proof}
Proposition~\ref{P:finite index for global fields}
yields $\widehat{A(K)} \isom \overline{A(K)}$,
so the first map is an isomorphism by Lemma~\ref{L:topological group}.
The second map is an isomorphism by \eqref{E:descent},
since $T\Sha$ is torsion-free by definition.
\end{proof}

The following is an analogue of Proposition~\ref{P:p-power}.

\begin{proposition}
\label{P:Selmer p-power}
If $A(K^s)[p^\infty]$ is finite,
then for any $v$, the map $\Selhat^{(p)} \to A(K_v)^{(p)}$
is injective.
\end{proposition}

\begin{proof}
Let $K_v' \subseteq K^s$ be the Henselization of $K$ at $v$.
Define ${\Sel^n}'$ and $\Selhat'$ in the same way as $\Sel^n$ and $\Selhat$,
but using $K_v'$ in place of its completion $K_v$.
By \cite{MilneADT}*{I.3.10(a)(ii)}, 
the natural maps ${\Sel^n}' \to \Sel^n$ and
$\Selhat' \to \Selhat$ are isomorphisms.
Similarly, by \cite{MilneADT}*{I.3.10(a)(i)} we may replace $A(K_v)^{(p)}$
by $A(K_v')^{(p)}:=\varprojlim A(K_v')/p^n A(K_v')$.
Hence it suffices to prove injectivity of 
$(\Selhat')^{(p)} \to A(K_v')^{(p)}$.

Choose $m$ such that $p^m A(K^s)[p^\infty] = 0$.
Suppose $b \in \ker\left( (\Selhat')^{(p)} \to A(K_v')^{(p)} \right)$.
For each $M \ge m$,
let $b_M$ be the image of $b$ in ${\Sel'}^{p^M} \subseteq H^1(K,A[p^M])$.
Then the image of $b_M$ under
\[
	{\Sel'}^{p^M} \to \frac{A(K_v')}{p^M A(K_v')} 
	\subseteq H^1(K_v',A[p^M]) \to H^1(K^s,A[p^M])
\]
is $0$.
The Hochschild-Serre spectral sequence 
(see \cite{MilneEtale}*{III.2.21})
gives an exact sequence
\[
	0 \to H^1(\Gal(K^s/K),A(K^s)[p^M]) \to H^1(K,A[p^M]) \to H^1(K^s,A[p^M])
\]
so $b_M$ comes from an element of $H^1(\Gal(K^s/K),A(K^s)[p^M])$,
which is killed by $\#A(K^s)[p^M] = p^m$.
Since $p^m b_M=0$ for all $M \ge m$,
we have $p^m b=0$.

By Lemma~\ref{L:Selmer torsion},
$b$ comes from a point in $A(K)[p^\infty]$,
which injects into $A(K_v')^{(p)}$.
\end{proof}

The next result is an analogue of Proposition~\ref{P:subscheme}.

\begin{proposition}
\label{P:Selmer subscheme}
Suppose that $A(K^s)[p^\infty]$ is finite.
Let $Z$ be a finite $K$-subscheme of $A$.
Then $Z(\A) \cap \Selhat = Z(K)$.
\end{proposition}

\begin{proof}
One inclusion is easy: $Z(K) \subseteq Z(\A)$
and $Z(K) \subseteq A(K) \subseteq \widehat{A(K)} \subseteq \Selhat$.
So we focus on the other inclusion.

As in the proof of \cite{Stoll2006preprint}*{Theorem~3.11}, 
we may replace $K$ by a finite extension 
to assume that $Z$ consists of a finite set of $K$-points of $A$.
Suppose $P \in Z(\A) \cap \Selhat$.
The $v$-adic component of $P$ in $Z(K_v)$
equals the image of a point $Q_v \in Z(K)$.
Then $P-Q_v$ maps to $0$ in $A(K_v)$, and in particular in $A(K_v)^{(p)}$,
so by Proposition~\ref{P:Selmer p-power} 
the image of $P-Q_v$ in $\Selhat^{(p)}$ is $0$.
This holds for every $v$, so if $v'$ is another place,
then $Q_{v'}-Q_v$ maps to $0$ in $\Selhat^{(p)}$.
The kernel of $A(K) \to \widehat{A(K)}^{(p)} \injects \Selhat^{(p)}$
is the prime-to-$p$ torsion of $A(K)$,
so $Q_{v'}-Q_v \in A(K)_\tors$.
This holds for all $v'$,
so the point $R:=P-Q_v$ belongs to $\Selhat_\tors$.
By Lemma~\ref{L:Selmer torsion}, $R \in A(K)_\tors$.
Thus $P=R+Q_v \in A(K)$.
So $P \in Z(\A) \intersect A(K) = Z(K)$.
\end{proof}

Each $\tau \in \Sel^{p^e}$ 
may be represented by a ``covering space'':
a locally trivial torsor $T$ under $A$
equipped with a morphism $\phi_T \colon T \to A$
that after base extension to $K^s$ 
becomes isomorphic to the base extension of $[p^e]\colon A \to A$.

\begin{lemma}
\label{L:p^e Selmer}
We have $\Selhat \subseteq \Union_{\tau \in \Sel^{p^e}} \phi_T(T(\A))$.
\end{lemma}

\begin{proof}
In the commutative diagram 
\[
\xymatrix{
& T(\A) \ar[d]^{\phi_T} \\
\Selhat \ar@{^{(}->}[r] \ar[d] & A(\A) \ar[d] \\
\Sel^{p^e} \ar@{^{(}->}[r] & \dfrac{A(\A)}{p^e A(\A)} \\
}
\]
the set $\phi_T(T(\A))$ is a coset of $p^e A(\A)$ in $A(\A)$,
and its image in $\frac{A(\A)}{p^e A(\A)}$ equals
the image of $\tau$ under the bottom horizontal map,
by the definition of this bottom map.
Thus any element of $\Selhat$ mapping to $\tau$ in $\Sel^{p^e}$
belongs to the corresponding $\phi_T(T(\A))$.
\end{proof}

\begin{proof}[Proof of Theorem~\ref{T:Selmer}]
We have
$X(K) \subseteq X(\A)^{\Br} \subseteq X(\A) \intersect \Selhat$,
by Proposition~\ref{P:BM for subvarieties}(a),
so it will suffice to show that $X(\A) \intersect \Selhat$
consists of $K$-rational points.

Let $e$ and $Y$ be as in Proposition~\ref{P:independent of F}.
For each $\tau \in \Sel^{p^e}$,
choose $\phi_T \colon T \to A$ as above,
choose $t \in T(K^s)$,
let $a = \phi_T(t) \in A(K^s)$,
and let $Y_a$ be the fiber of $Y \to A$ above $a$, viewed as
a finite subscheme of $X_{K^s}$.
Choose a finite $K$-subscheme $Z$ of $X$ such that $Z_{K^s}$ contains
all these $Y_a$ (as $\tau$ ranges over $\Sel^{p^e}$).

For each $v$, choose a separable closure $K_v^s$ of $K_v$ containing $K^s$.
By the proof of Lemma~\ref{L:exists Z},
any $p$-basis of $K$ is also a $p$-basis for $K_v$,
and hence for $K_v^s$,
so the conclusion of Proposition~\ref{P:independent of F} may be applied
with $F=K_v^s$ to yield
\[
	X(K_v^s) \intersect (a + p^e A(K_v^s)) \subseteq Y_a(K_v^s)
\]
for each $a$.
By definition of $\phi_T$, we have $\phi_T(T(K_v^s))=a + p^e A(K_v^s)$.
Thus
\[
	X(K_v) \intersect \phi_T(T(K_v)) \subseteq X(K_v^s) \intersect (a + p^e A(K_v^s)) \subseteq Y_a(K_v^s) \subseteq Z(K_v^s).
\]
Hence
\[
	X(K_v) \intersect \phi_T(T(K_v)) \subseteq X(K_v) \intersect Z(K_v^s) = Z(K_v).
\]
This holds for all $v$, so 
$X(\A) \intersect \phi_T(T(\A)) \subseteq Z(\A)$.
Taking the union over $\tau \in \Sel^{p^e}$, 
and applying Lemma~\ref{L:p^e Selmer},
we obtain
$X(\A) \intersect \Selhat \subseteq Z(\A)$.
Thus
$X(\A) \intersect \Selhat \subseteq Z(\A) \intersect \Selhat$,
and the latter equals $Z(K)$ by Proposition~\ref{P:Selmer subscheme},
so we are done.
\end{proof}

\section*{Acknowledgements}
This work was initiated while the authors were 
at the program on ``Rational and integral points on
higher-dimensional varieties'' at the
Mathematical Sciences Research Institute in Spring 2006.
We thank 
Zo\'e Chatzidakis, 
Jean-Louis Colliot-Th\'el\`ene, 
Fran\c{c}oise Delon, 
Cristian Gonzalez-Aviles, 
Tom Scanlon,
and 
Michael Stoll 
for discussions.
B.P.\ was partially supported by NSF grant DMS-0301280.

\end{document}